\documentclass[fontsize=12pt,a4paper,headings=normal,
twoside=false,leqnoparskip=half-,abstract=true]{scrartcl}
\usepackage[english]{babel}
\usepackage[utf8]{inputenc}
\usepackage{hyperref}
\usepackage{blkarray}
\usepackage{comment}
\hypersetup{
 pdftitle={},
 pdfauthor={},
 colorlinks=true,
 linkcolor=blue,
 citecolor=blue,
 filecolor=blue,
 urlcolor=blue}

\usepackage[pagewise]{lineno}%\linenumbers
\usepackage[version=4]{mhchem}
\usepackage[format=plain,labelfont=bf,font=small]{caption}
\usepackage{subfigure}%replaced by subcaption
\usepackage{xcolor}
\usepackage[arrow, matrix, curve]{xy}
\usepackage{tikz}
\usepackage{float}
\usepackage{orcidlink}

\usepackage{tikz-cd}
\usepackage{caption}
\usepackage{amsmath,amsthm}
\usepackage{amssymb} 
\usepackage[normalem]{ulem}
\usepackage{academicons}
\definecolor{orcidlogocol}{HTML}{A6CE39}

\newtheorem{theorem}{Theorem}[section]
\newtheorem{lemma}[theorem]{Lemma}
\newtheorem{corollary}[theorem]{Corollary}
\newtheorem{proposition}[theorem]{Proposition}

\theoremstyle{definition}
\newtheorem{definition}[theorem]{Definition}

\newcommand{\be}{\begin{equation}}
\newcommand{\ee}{\end{equation}}
\newcommand{\bea}{\begin{eqnarray}}
\newcommand{\eea}{\end{eqnarray}}

\theoremstyle{remark}
\newtheorem{remark}[theorem]{Remark}

\numberwithin{equation}{section}

\DeclareMathAlphabet{\mathpzc}{OT1}{pzc}{m}{it}

\date{\today}

\graphicspath{{./images/} }

\begin{document}

\title{Mass action systems: two criteria for Hopf bifurcation without Hurwitz}
%\author{Alex Blokhuis \orcidlink{0000-0002-4594-596X}}
%\email{alex.blokhuis@hotmail.com}
%\affiliation{University of Strasbourg $\&$ CNRS, UMR7140, 67083 Strasbourg, France} 
%\affiliation{Instituto IMDEA Nanociencia, Calle Faraday 9, 28049 Madrid, Spain} 
\author{Nicola Vassena \orcidlink{0000-0001-5411-4976}}
%\affiliation{Interdisciplinary Center for Bioinformatics, Leipzig University, Germany, nicola.vassena@uni-leipzig.de} 

\maketitle

\begin{abstract}
We state two sufficient criteria for periodic oscillations in mass action systems. Neither criterion requires a computation of the Hurwitz determinants. Instead, both criteria exploit the linear algebra concepts of $D$-stability and $P$-matrices. The criteria are complementary: the first is based on a stable matrix that is not a $P^-$ matrix, while the second is based on a $P^-$ matrix that is not stable. In analogy, a qualitatively different interpretation follows: the first criterion relates to positive feedback in the network, while the second concerns negative feedback. We present examples that showcase the applicability of both criteria. As a final independent remark, we prove that for the special case of fully-open networks, the capacity for Hopf bifurcation is just equivalent to the capacity for a steady-state with a complex pair of eigenvalues with positive-real part.

\end{abstract}

\tableofcontents

\section{Introduction}\label{sec:intro}
The quest of finding periodic orbits in polynomial differential equations is notoriously difficult. A standard approach aims at proving the occurrence of a Hopf bifurcation, where an equilibrium changes stability and generates a periodic orbit. Purely imaginary eigenvalues of the Jacobian are a necessary spectral condition for a Hopf bifurcation to occur. Therefore, the great majority of the literature in reaction networks employs an explicit Hurwitz criterion \cite{Liu94}. See among others \cite{Gat2005, DomKirk09, Hell16, HMNPSHopf17, conradietal19, TangHopf22, BaBo23}. The major drawback of the Hurwitz approach is its computational complexity, which limits its applicability to small networks. Moreover, such a computation may obfuscate, rather than illuminate, the underlying chemical mechanism that leads to oscillations.

With these motivations, we present two criteria that establish nonstationary periodic oscillations in mass action systems and do not require any Hurwitz computation. The first ingredient is the standard Jacobian parametrization by Clarke \cite{ClarkeSNA}, which expresses the Jacobian matrix $Jac$ at steady states in the form of
\begin{equation}
Jac=AD.
\end{equation}
{In Clarke's work, the matrix \( A \) is reparametrized using the extreme ray of the convex cone of steady-state fluxes. For this reason, this reparametrization at steady states is often referred to as \emph{convex coordinates}. More relevant here, the matrix \( D \) is a parametric positive diagonal matrix, with parameters related to the steady-state concentrations. It is then possible to analyze the stability of the Jacobian \( AD \) by varying the parameters in \( D \). {In the case of $A$ being invertible, 
\begin{equation}
\operatorname{rank} AD = \operatorname{rank} A
\end{equation}}
ensures {that $AD$ is invertible for all choices of \( D \)}. Therefore, any change in stability must involve a crossing at purely imaginary eigenvalues. The second ingredient are the linear algebra concepts of \( D \)-stability and \( P^- \) matrices: \( D \)-stability generalizes matrix stability, requiring that a matrix \( A \) remains stable under multiplication by any positive diagonal matrix \( D \). \( P^- \) matrices are those whose principal minors follow a prescribed sign pattern: all \( k \)-principal minors have sign \( (-1)^k \). In particular, we use two complementary results: first, the straightforward observation that any matrix with a \( k \)-principal minor of sign \( (-1)^{k-1} \) is not \( D \)-stable; and second, a classic result by Fisher and Fuller \cite{FisherFuller58, Fisher72simple}, which guarantees the existence of a positive diagonal matrix \( D \) such that \( AD \) is stable whenever \( A \) is a \( P^- \) matrix.} We rely on such results to provide manageable sufficient conditions on $A$ for the existence of a positive diagonal matrix $D^*$ such that $AD^*$ has purely imaginary eigenvalues. Finally, we return to nonlinear dynamics: the third ingredient is a result by Fiedler \cite[Theorem 4.7]{Fiedler85PhD} on \emph{global Hopf bifurcation} that - together with the first two ingredients - guarantees the existence of nonstationary periodic orbits.

The paper is organized as follows. Section \ref{sec:rn} introduces setting and notation. Section \ref{sec:Jacobian} reviews the parametrization method Stoichiometric Network Analysis by Clarke. Section \ref{sec:linearalgebra} is dedicated to linear algebra and presents in more detail the concepts of $D$-stability, $P^-$-matrices, and the Fisher-Fuller result. Section \ref{sec:main} contains the main Theorem \ref{thm:2} and - as Corollary \ref{cor:criteria} - the two criteria. Section \ref{sec:interpretation} provides an interpretation of the chemical mechanism underlying the criteria: Criterion I concerns an unstable-positive feedback within a stable network; Criterion II concerns an unstable-negative feedback. We present examples for both criteria in Section \ref{sec:examples}. Finally, Section \ref{sec:fullyopen} provides an independent result for the special case of fully-open networks. Section \ref{sec:conclusion} summarizes the paper and outlooks future directions.

%The first criterion is based on \emph{global} Hopf bifurcation, as developed originally by Yorke et al. \cite{AlexYorke78, ChowMPYorkeGH78, MpYorke82, AllYorke84}, and extended by Fiedler \cite{Fiedler85PhD, F19}, and relies on the concepts of $D$-stability and $P$-matrices. The second criterion is based on \emph{local} Hopf bifurcation, see e.g. \cite{GuHo84}.

%A recent body of \emph{inheritance} results by Banaji and co-authors \cite{ba} aims at lifting properties . For large-size networks, it is not yet available nor foreseeable an efficient algorithm that checks on the presence of reduced
%On the other hand, 
%A key-observation is that the source of interest in the  GLOBAL HOPF. 

%Reaction networks are often open systems, exchanging chemicals with the outside environment via inflow and outflow reactions. Several results in the literature have focused on fully-open systems: for instance.

\section{Reaction networks}\label{sec:rn}

A reaction network $\mathbf{\Gamma}=(S,R)$ consists of a set $S$ of species that interact via a set $R$ of reactions. A reaction $i$ is an ordered association between \emph{reactant} and \emph{product} species:
\begin{equation} \label{reactionj}
 i: \quad s^{i}_1X_1+...+s^{i}_{|S|}X_{|S|} \underset{i}{\longrightarrow} \tilde{s}^{i}_1X_1+...+\tilde{s}^{i}_{|S|}X_{|S|}.
\end{equation}
Here, $X_1, ..., X_{|S|}$ indicate $|S|$ distinct species and the nonnegative integer coefficients $s^{i}_m,\tilde{s}^i_m$  are the so-called \emph{stoichiometric coefficients}. A reaction without reactants and single product $X_m$ is called \emph{an inflow reaction} to $X_m$. Respectively, a reaction without products and a single reactant $X_m$ is called \emph{an outflow reaction} from $X_m$.
{All reactions that are neither inflow nor outflow reactions are assumed to have both reactants and products.} A network is called \emph{closed} if there are no inflow and outflow reactions. At the other extreme, a network whose reaction set $R$ contains inflow reactions to all its species and outflow reactions from all of its species is called \emph{fully open}. Finally, a reaction for which a species $X_m$ is both a reactant and a product is called \emph{explicitly autocatalytic}.

Consider the $|S|$-vector $x>0$ of the positive species concentrations in a well-mixed reactor. The dynamics $x(t)$  follows the system of Ordinary Differential Equations:
\begin{equation}\label{maineq}
    \dot{x}=f(x):=N \mathbf{r}(x),
\end{equation}
where $N$ is the $|S|\times |R|$ \emph{stoichiometric matrix}, defined by
\begin{equation}\label{S}
    N_{mi}:=\tilde{s}^i_m - s^i_m,
\end{equation} and $\mathbf{r}(x)$ is the $|R|$-vector of \emph{reaction rates} (kinetics). Equilibria vector $\bar{x}$ with $$0=f(\bar{x})$$ are called \emph{steady-states} of $\mathbf{\Gamma}$. As our criteria are based on the bifurcation of steady-states, we consider throughout only networks $\mathbf{\Gamma}$ that admit at least one positive steady-state. This condition just requires the existence of at least one positive vector $\mathbf{v}\in \mathbb{R}^{|R|}_{>0}$, called here \emph{steady-state flux vector} such that
\begin{equation}\label{eq:fluxvector}
N\mathbf{v}=0,
\end{equation}
i.e., $\mathbf{v}$ is a positive right kernel vector of the stoichiometric matrix. Networks, whose stoichiometric matrix satisfies \eqref{eq:fluxvector} for a positive vector $\mathbf{v}>0$, have been called \emph{consistent} \cite{Ang07} or also \emph{dynamically nontrivial} \cite{BaBo23}.

It is reasonable to assume that each reaction rate $r_i$ is a positive monotone function of the concentrations $x_m>0$ of those species $X_m$, which are reactants to $i$, i.e., for $s^i_m>0$. This class of nonlinearities has been named \emph{monotone chemical functions} \cite{VasHunt} or \emph{weakly monotone kinetics} \cite{ShiFei12} in the literature. A most prominent example is mass action kinetics \cite{MA64}, which considers any reaction rate $r_i$ as a monomial function of the concentrations:
\begin{equation}
  r_i(x) := a_i\prod_{m=1}^{|S|} x_m^{s_{m}^{i}}. 
\end{equation}
The coefficient $a_i>0$ is a positive parameter, and the exponents $s_{m}^{i}$ are the stoichiometric coefficients of species $X_m$ as a reactant of the reaction $i$ in \eqref{reactionj}. Note that the rate of an inflow reaction is then constant, i.e., $r_i(x)=a_i$. More general polynomial functions such as Generalized Mass Action kinetics \cite{Muller:12} or rational functions such as  Michaelis-Menten \cite{MM13} and Hill \cite{Hill10} kinetics also fall within the class of monotone chemical functions.

%It has been already observed the deceptive nature of mass action kinetics in being the `simplest case' of kinetics. On one hand, few parameters and more established results and tools from algebraic geometry and computer algebra suggest mass action as the first simple case. On the other hand, the concept of `parameter-rich' kinetics, which excludes mass action, has been recently introduced \cite{VasStad23} to underline how the questions of stability and bifurcations are more easily treatable with more parameters.

\section{Jacobian parametrization}\label{sec:Jacobian}
Bruce L. Clarke developed \emph{Stoichiometric Network Analysis} in a series of papers: see \cite{ClarkeSNA} and references therein. In the very same paper \cite{ClarkeSNA}, Clarke also commented on one possibility to circumvent the computation of Hurwitz determinants by studying the $D$-stability of the Jacobian matrix in convex coordinates. This comment highly resonates with the content of our main Theorem \ref{thm:2} below. Clarke further cited a paper of his with some announced related network results, submitted to `\emph{Linear Algebra and Applications}'. To my best effort, however, I have been unable to retrieve such a paper and even to confirm its publication anywhere: it is not listed as a paper in the database of \emph{Linear Algebra and its Applications} nor in Clarke's publication list. It is nevertheless possible that what we state here connects to what Clarke had in mind. 

We introduce only the necessary standard technicalities, we again refer to \cite{ClarkeSNA} for a thorough overview of the subject.
Consider the \emph{steady-state flux cone} $\mathpzc{F}$, i.e. the set of steady-state flux vectors:
$$\mathpzc{F}\quad:=\quad\{\mathbf{v}\in \mathbb{R}_{>0}^{|R|} \quad | \quad N\mathbf{v}=0\}.$$
 In his work, Clarke uses convex coordinates to parametrize $\mathpzc{F}$. In particular, via convexity, we can consider the extreme rays of $\mathpzc{F}$, $\{E_1,...,E_p\}$, and express each steady-state flux $\mathbf{v}\in\mathpzc{F}$ as a {nonnegative} linear combination of such extreme rays: 
$$\mathbf{v}=E \mathbf{j},$$
where $\mathbf{j}\in \mathbb{R}^p_{\ge0} $, and the $|R|\times p$ \emph{extreme-ray matrix} $E$ has columns $\{E_1,...,E_p\}$. The chosen parametrization of the flux cone is anyway not relevant in the following arguments.

The relevant observation is perhaps only that the derivative of a univariate polynomial $p(x)=ax^n$ can be expressed as
\begin{equation}\label{eq:polder}
p'(x)=n\;ax^{n-1}=n\;ax^{n-1}\;\frac{x}{x}=n p(x) \frac{1}{x}.
\end{equation}
This implies that the Jacobian matrix $Jac$ of mass action systems \eqref{maineq} at a steady-state $\bar{x}$ reads:
\begin{equation}
\begin{split}
Jac=f_x(x)|_{x=\bar{x}}=N\frac{\mathbf{r}(x)}{\partial x}\bigg|_{x=\bar{x}}=
N \frac{\mathbf{r}(x)}{\partial x}\bigg|_{x=\bar{x}} \operatorname{diag}\bigg(\frac{\bar{x}}{\bar{x}}\bigg)
= B(\bar{\mathbf{v}})\operatorname{diag}\bigg(\frac{1}{\bar{x}}\bigg),
\end{split}
\end{equation}
where $B(\bar{\mathbf{v}})$ only depends on the choice of $\bar{\mathbf{v}}$ and not on the value of $\bar{x}$ anymore. Following Clarke, a straightforward computation generalizes \eqref{eq:polder} and shows that 
\begin{equation}
B(\bar{\mathbf{v}})=N\operatorname{diag}(\bar{\mathbf{v}})Y^T,
\end{equation}
where the $|S|\times |R|$ \emph{kinetic matrix} $Y$ is defined as $$Y_{mi}:=s_{i}^m.$$ $Y$ encodes the exponents of the monomials in \eqref{maineq} endowed with mass action kinetics.

Note that the values $1/\bar{x}_i$ can be thought of as positive parameters themselves, $h_i:=1/\bar{x}_i$. Thus, the theory of Stoichiometric Network Analysis identifies two sets of parameters: $( \mathbf{h},\mathbf{j})$.
One parameter set, $\mathbf{j}$, determines the steady-state flux vector $\bar{\mathbf{v}}$; the other parameter set, $\mathbf{h}$, determines the steady-state concentrations $\bar{x}$. {The parameters $(\mathbf{h},\mathbf{j})$ are called \emph{steady-state parameters}}. The Jacobian matrix can then be expressed {in steady-state parameters} as
$$Jac(\mathbf{h},\mathbf{j})=N\operatorname{diag}E\mathbf{j}\operatorname{diag}(\mathbf{h})$$
or, having chosen $\bar{\mathbf{v}}$ a priori, simply
$$Jac(\mathbf{h},\bar{\mathbf{v}})=B(\bar{\mathbf{v}})\operatorname{diag}(\mathbf{h}).$$

{It is straightforward to recover the reaction rate constants and - if present - the values of the conserved quantities from any fixed choice of steady-state parameters $(\bar{\mathbf{h}},\bar{\mathbf{j}})$, we refer again to Clarke \cite{ClarkeSNA} for explicit conversion formulas. In particular, the existence of steady-state parameters for which a certain dynamical feature occurs directly implies the existence of reaction rate constants for which the same dynamical feature occurs.}

\section{Linear algebra: $D$-stability and $P^-$ matrices}\label{sec:linearalgebra}
This section reviews the standard concepts of $D$-stability and $P^-$ matrices. We refer to \cite{Gio15, Ku19} for a more extended overview. We start with the definition of \emph{inertia} of a matrix.
\begin{definition}[Inertia of a matrix]
    The inertia of an $n \times n$ square matrix $A$ is a nonnegative triple  $$\operatorname{inertia}(A):=(\sigma^-_A,\sigma^+_A,\sigma^0_A),$$
where $\sigma^-_A$,$\sigma^+_A$ and $\sigma^0_A$
 are the number of eigenvalues of $A$ with negative-real part, positive-real part, and zero-real-part eigenvalues, respectively. The eigenvalues are counted with their multiplicities so that $\sigma^+_A+\sigma^-_A+\sigma^0_A=n$.
 \end{definition}

\begin{definition}[(in)stability, $D$-(in)stability]
A $n\times n$ matrix $A$ is \emph{stable} if all of its eigenvalues have negative real part.
Conversely, $A$ is \emph{unstable} if at least one eigenvalue has a positive real part. Moreover, a matrix $A$ is \emph{$D$-stable} if $AD$ is stable for any choice of a positive diagonal matrix $D$. Conversely, $A$ is  \emph{$D$-unstable} if there exists a diagonal matrix $D$ such that $AD$ is unstable.
\end{definition} 

{Clearly, the inertia of a stable matrix is  
$(n,0,0)$. Note also that a matrix can be neither ($D$-)stable nor ($D$-)unstable, e.g. the zero matrix. More importantly,} by choosing $D=\operatorname{Id}$, it is clear that $D$-stability implies stability, and that instability implies $D$-instability. In addition, the $D$-stability of $A$ necessarily requires the $D$-stability of all its principal submatrices, as the next Lemma shows. We fix the notation: let $\kappa$ be any choice of $k\le n$ indices in $\{1,...,n\}$. $A[\kappa]$ denotes the principal submatrix of $A$ with column/row index $\kappa$; {its determinant is called \emph{a $k$-principal minor} of $A$.}
\begin{lemma}\label{lem:Dinstabilitynec}
If any of the principal submatrices of $A$ is $D$-unstable, then $A$ is $D$-unstable.
\end{lemma}
\proof
 By assumption, there exists a choice of $\kappa$ such that $A[\kappa]$ is $D$-unstable. Without loss of generality, we can assume $\kappa=\{1,...,k\}$. Choose now $D[\kappa]=\operatorname{diag}(d_1,...,d_k)$ such that $A[\kappa]D[\kappa]$ is unstable, and rescale all other entries of $D$, $d_i$ with $i>k$, as $d_i=\varepsilon$. Consider now the family of matrices $AD(\varepsilon)$. Clearly, for $\varepsilon=0$, the spectrum of $AD(0)$ corresponds to the spectrum of $A[\kappa]D[\kappa]$ plus $n-k$ eigenvalues zero. In particular $AD(0)$ is unstable by construction. The continuity of the eigenvalues with respect to the entries guarantees that such instability persists for $\varepsilon$ small enough, which yields the instability of $AD$, for a choice of positive $D$, and thus the $D$-instability of $A$.
\endproof

Unfortunately, proving that a square matrix $A$ is $D$-stable is computationally nontrivial. It may precisely reduce to exclude purely-imaginary eigenvalues of $AD$ for any choice of $D$, via the Hurwitz computation. To avoid being a dog chasing its tail, we have to rely on sufficient conditions on $A$ that guarantee the existence of positive diagonal matrices $D_1$ and $D_2$ such that 
\begin{equation}\label{eq:diffinertia}
    \operatorname{inertia}(AD_1)\neq\operatorname{inertia}(AD_2),
\end{equation}
which sufficiently guarantees that $A$ is not $D$-stable. The key here is the concept of $P^-$ matrices. 
\begin{definition}[$P^-$ and $P^-_0$ matrices]\label{def:Pmatrix}
A $n\times n$ matrix $A$ is called a $P^-$ matrix if all of its $k$-principal minors have sign $(-1)^k$. A $n\times n$ matrix $A$ is called a $P^-_0$ matrix if all of its \emph{nonzero} $k$-principal minors have sign $(-1)^k$.
\end{definition}
The set of $P^-_0$ matrices is just the closure of the open set of $P^-$ matrices. We recall two standard results that relate $D$-stability and $P^-_0$ matrices. The first proposition follows directly from Lemma \ref{lem:Dinstabilitynec}.
\begin{proposition}\label{prop:Pf}
Let $A$ be a $n\times n$ matrix that is not a $P^-_0$ matrix. Then there exists a positive diagonal $D$ such that $AD$ is unstable. 
\end{proposition}
\proof
If $A$ is not a $P^-_0$ matrix, then $A$ possesses a principal submatrix $A[\kappa]$ such that { $$\operatorname{sign}\operatorname{det}A[\kappa]=(-1)^{k-1}.$$}
In particular $A[\kappa]$ possesses an odd number of real positive eigenvalues, and it is thus unstable. Recalling that instability implies $D$-instability, Lemma \ref{lem:Dinstabilitynec} yields the statement.
\endproof

The second is a celebrated result by Michael E. Fisher and A. T. Fuller in 1958 \cite{FisherFuller58}, elaborated also by Franklin M. Fisher in \cite{Fisher72simple}. For the full generality of this theorem, a further definition is required.
\begin{definition}[Fisher and Fuller $P^-_{FF}$ matrices]
A $P^-_0$ matrix $A$ is called a $P^-_{FF}$ matrix, or \emph{Fisher\&Fuller matrix}, if there is at least a sequence of nested nonsingular principal matrices $$(A[\kappa_1],A[\kappa_2],...,A[\kappa_{n-1}],A),$$
of every order $|\kappa_i|=i=(1,...,n)$, such that $A[\kappa_{i-1}]$ is a principal submatrix of $A[\kappa_{i}]$.
\end{definition}
Clearly, the following set inclusion holds:
$$\text{$P^-$ matrices}\quad\subset\quad\text{$P^-_{FF}$ matrices}\quad\subset\quad\text{$P^-_{0}$ matrices}.$$
A warning: Fisher \cite{Fisher72simple} uses the word \emph{Hicksian} to refer to $P^-_{FF}$ matrices. Such a name can be also found in the literature \cite{Gio15} to refer to $P^-$ matrices.
\begin{theorem}[Theorem 1' in \cite{Fisher72simple}]\label{thm:FF}
Let $A$ be a $P^-_{FF}$ matrix. Then there exists a positive diagonal $D$ such that $AD$ has all eigenvalues that are real, negative, and simple.
\end{theorem}
We can derive easy conditions for \eqref{eq:diffinertia} based on Proposition \ref{prop:Pf} and Theorem \ref{thm:FF}. We do this explicitly in the next two corollaries.
\begin{corollary}[of Proposition \ref{prop:Pf}]\label{cor:1}
    Consider a $n\times n$ matrix $A$ that is not a $P^-_0$ matrix. Further, assume that $A$ has no eigenvalue with a positive-real part. Then there exist two diagonal matrices $D_1$ and $D_2$ such that \eqref{eq:diffinertia} holds.
\end{corollary}
\proof
{Consider $D_1=\operatorname{Id}$ and let $D_2$ be the matrix satisfying the statement of Proposition \ref{prop:Pf}. The statement follows from noting that
$AD_1=A\operatorname{Id}=A$ has no eigenvalues with positive real part while $AD_2$ has an eigenvalue with positive real part via Proposition \ref{prop:Pf}. 
\endproof}
\begin{corollary}[of Theorem \ref{thm:FF}]\label{cor:2}
Consider an unstable $P^-$ matrix $A$. Then there exist two diagonal matrices $D_1$ and $D_2$ such that \eqref{eq:diffinertia} holds.
\end{corollary}
\proof
{Consider again $D_1=\operatorname{Id}$ and let $D_2$ be the matrix satisfying the statement of Theorem \ref{thm:FF}.
The statement follows by noting that $AD_1=A\operatorname{Id}=A$ has an eigenvalue with positive real part while $AD_2$ has only eigenvalues with negative real part via Theorem \ref{thm:FF}.
\endproof}

\paragraph{Purely imaginary eigenvalues of $AD$.} {In the case in which $A$ is invertible}, the condition \eqref{eq:diffinertia} implies the existence of a positive matrix $D^*$ such that $AD^*$ has purely imaginary eigenvalues. We formally argue as follows. Let $\mathbf{d}_1$ and $\mathbf{d}_2$ indicate the vectors in $\mathbb{R}^n_{>0}$ such that $D_1=\operatorname{diag}(\mathbf{d}_1)$ and $D_2=\operatorname{diag}(\mathbf{d}_2)$. Consider any regular curve 
$$\gamma(t):[0,1] \mapsto {\mathbb{R}^n},$$
such that $\gamma(0)=\mathbf{d}_1$ and $\gamma(1)=\mathbf{d}_2$, and the associated \emph{eigenvalue curve}
$$\Lambda(t):[0,1]\mapsto {\mathbb{R}^n}\quad \text{ defined by}\quad \Lambda(t)=(\lambda_1(t),...,\lambda_n(t))=\operatorname{eigenvalues}(A\operatorname{diag}(\gamma(t)).$$
The condition \eqref{eq:diffinertia} implies a change in the sign of the real part of at least one of the eigenvalues along the path $\gamma$: intermediate value theorem guarantees then the existence of at least one $t^*$ such that $\Lambda(t^*)$ identifies at least one eigenvalue with zero-real part. Note that {
$$\operatorname{rank}A\operatorname{diag}(\gamma(t))=\operatorname{rank}A$$}
prevents real zero eigenvalues, and thus $t^*$ identifies at least one pair of purely imaginary eigenvalues. We conclude this section with an open linear algebra question:\\ 

\begin{center}
$\mathbf{Q^*}:$  Does it always exist a choice of path $\gamma$ such that \\
    (i) $A\operatorname{diag}(\gamma(t^*))$ has a \emph{simple} pair of purely imaginary eigenvalues $\lambda_{i,j}$,\\ 
    (ii) $\Re(\lambda_{i,j}(t^*))'\neq 0$?
\end{center}
\textcolor{white}{h}

{In Section \ref{sec:Jacobian}, we have shown} how a Jacobian matrix evaluated at steady-states of mass action systems can be expressed precisely as a product $AD$. A positive answer to $\mathbf{Q^*}$ would guarantee the applicability of the standard local Hopf theorem \cite{GuHo84} (or Theorem \ref{thm:Hopf} in Section \ref{sec:fullyopen}) to prove the existence of periodic orbits. As we do not know the answer to $\mathbf{Q^*}$ yet, we will argue instead via \emph{global Hopf bifurcation}. 
We proceed to our main result.

\section{Nonlinear dynamics: one theorem and two criteria}\label{sec:main}

 This section glues the results from Sections  \ref{sec:Jacobian} and \ref{sec:linearalgebra} to establish periodic orbits for mass action systems without any Hurwitz computation. The central nonlinear argument relies on the theory of global Hopf bifurcation. This body of work appears to be essentially unknown in reaction network literature, which is mostly concerned with local methods based on Hurwitz computation. A relevant exception is the work \cite{F19} by Fiedler, which however focuses on the restricted stoichiometric structure of feedback cycles and does not cover mass action kinetics.

A short guide for neophytes: Alexander \& Yorke \cite{AlexYorke78} addressed the possibility of global bifurcation of periodic orbits. In particular, they established periodic orbits whenever an \emph{odd} number of pairs of complex conjugate eigenvalues crosses the imaginary axis, at an equilibrium with an invertible Jacobian. A joint effort by Chow, Mallet-Paret, and Yorke \cite{ChowMPYorkeGH78} extended the very same result to include any net change of stability. The description of the global continua of periodic orbits, so-called \emph{snakes}, have been addressed in \cite{MpYorke82} for a generic situation. Alligood and Yorke \cite{AllYorke84} lifted the result in the non-generic case. Via a further analyticity assumption, Fiedler's work \cite{Fiedler85PhD} on general parabolic systems included Hopf points that are not necessarily isolated: these are the generalities we consider in the following Theorem \ref{thm:2}. {Let us first recall Fiedler's result in a version that is consistent with our setting.
\begin{theorem}[(4.7) from \cite{Fiedler85PhD}]\label{thm:fiedler} Consider a parametric vector field depending on a scalar parameter $\lambda$: 
$$\dot{x}=f(x,\lambda),$$
with $x \in \mathbb{R}^N, \lambda\in [a,b] \subseteq \mathbb{R}$. Assume the following conditions:
\begin{enumerate}
\item Analyticity of $f$ in $x$ and $\lambda$;
\item A family of steady-states is parametrized over $\lambda$, i.e. there exists $\bar{x}(\lambda)$ s.t.
$$0=f(\bar{x}(\lambda),\lambda)\quad\text{for}\quad \lambda \in [a,b];$$
\item The Jacobian of $\bar{x}(\lambda)$ is invertible for any $\lambda$:
$$\operatorname{det}f_x(x,\lambda)|_{(\bar{x}(\lambda),\lambda)}\neq 0 \quad \text{for}\quad \lambda\in[a,b];$$
\item Net-change of stability of the Jacobian:
$$\operatorname{inertia}f_x(x,\lambda)|_{(\bar{x}(a),a)}\neq \operatorname{inertia}f_x(x,\lambda)|_{(\bar{x}(b),b)}.$$
\end{enumerate}
Then there exists $\lambda^*$ such that $\dot{x}=f(x,\lambda^*)$  admits nonstationary periodic orbits.
\end{theorem}
We note that \cite{Fiedler85PhD} presents Theorem \ref{thm:fiedler} in a broader setting for parabolic systems and includes stronger conclusions related to the theory of global Hopf bifurcations, which we have omitted here for simplicity of presentation.} We state now our main result.

\begin{theorem}\label{thm:2}
Let $\mathbf{\Gamma}$ be a network with stoichiometric matrix $N$ and kinetic matrix $Y$. Assume $N$ has full-rank $|S|$ and consider a steady-state flux vector $\bar{\mathbf{v}}>0$ such that  
$$A:=B(\bar{\mathbf{v}})=N\operatorname{diag}(\bar{\mathbf{v}}) Y^T,$$ is invertible. Assume moreover that there exist two positive diagonal matrices $D_1$ and $D_2$ such that \eqref{eq:diffinertia} holds, i.e.,
$$\operatorname{inertia}AD_1\neq \operatorname{inertia}AD_2.$$
Then there exists a choice of reaction rates {constants} such that the associated mass action system \eqref{maineq} admits nonstationary periodic solutions. 
\end{theorem}

\proof
Recall the parameters $h_i=1/\bar{x}_1$. 
We have
$$\operatorname{det}(Jac( \mathbf{h},\bar{\mathbf{v}}))=\operatorname{det}(B(\bar{\mathbf{v}})\operatorname{diag}(\mathbf{h}))=\operatorname{det}(B(\bar{\mathbf{v}}))\operatorname{det}(\operatorname{diag}(\mathbf{h}))\neq 0,$$
that is, the Jacobian $Jac( \mathbf{h},\bar{\mathbf{v}})$ is invertible for any choice of $\mathbf{h}$.

Consider now the two choices of $\textbf{h}_i$, $i=1,2$ such that
$$\operatorname{diag}(\textbf{h}_1)=D_1\quad\quad  \text{and}\quad\quad \operatorname{diag}(\textbf{h}_2)=D_2.$$
In analogy to what we discussed in Section \ref{sec:linearalgebra}, we may consider any \emph{analytic} curve $\gamma(\beta)$ in $\mathbb{R}^{|S|}_{>0}$, $\beta\in[0,1]$ with $\gamma(0)=\mathbf{h}_1$, $\gamma(1)=\mathbf{h}_2$. That is, $\gamma$ connects $\mathbf{h}_1$ to $\mathbf{h}_2$ and thus parametrizes a family of steady-states of \eqref{maineq} with invertible Jacobian and a change of inertia from $\beta=0$ to $\beta=1$. Due to the analyticity of mass action systems and the curve $\gamma(\mathbf{h})$, we can apply Theorem {\ref{thm:fiedler}} and global Hopf bifurcation yields nonstationary periodic orbits. 
\endproof

Based on the abstract Theorem \ref{thm:2}, we simply use Corollaries \ref{cor:1} and \ref{cor:2} to state computationally manageable sufficient criteria for periodic orbits in mass action systems. Criterion I builds on Corollary \ref{cor:1}, and Criterion II on Corollary \ref{cor:2}.

\begin{corollary}\label{cor:criteria}
Assume there exists a steady-state flux vector $\bar{\mathbf{v}}$ such that one of the two following conditions holds.
\begin{enumerate}
\item \textbf{Criterion I:} $B(\bar{\mathbf{v}})$ is stable and not a $P^-_0$ matrix;
\item \textbf{Criterion II:} $B(\bar{\mathbf{v}})$ is an unstable $P^-_{FF}$ matrix.
\end{enumerate}
then there exists a choice of reaction rates {constants} such that the associated mass action system \eqref{maineq} admits nonstationary periodic solutions.
\end{corollary}

\begin{remark}[Conserved quantities] 
For simplicity of presentation, Theorem \ref{thm:2} and the two sufficient criteria verbatim apply only to systems without conserved quantities. In fact, the presence of conserved quantities prevents $B(\bar{\mathbf{v}})$ from being invertible, and thus in our definitions $B(\bar{\mathbf{v}})$ cannot be stable nor a $P^-_{FF}$ matrix. It is clear that this is no issue, as the results can be identically applied to any invertible principal submatrix of $B(\bar{\mathbf{v}})$, via a straightforward perturbation argument. 
%See also  for \textcolor{blue}{the full-rank} \emph{reduced Jacobian} in each stoichiometric compatibility class, after a standard reduction procedure \cite{BaPa16,CoCa19}. 
\end{remark}
\begin{remark}
$P^{-}_0$ matrices as Jacobians have been addressed in reaction networks in \cite{Ba-07, VasHunt}. In the first contribution \cite{Ba-07}, it is shown that a $P^-_0$ Jacobian implies the \emph{injective property} and thus excludes multistationarity, i.e., the coexistence of multiple steady-states in the same stoichiometric compatibility class. Criterion II above, in particular, describes a mechanism for periodic oscillations in monostationary networks. See \cite{toth99} for a discussion with examples of the relationship between periodic oscillations and multistationarity in chemical reaction networks.
 \end{remark}

\section{Interpretation: Positive and negative feedbacks}\label{sec:interpretation}
In biochemical systems, periodic oscillations have long been associated with the presence of \emph{positive and negative feedbacks}, see \cite{tyson1975} for a review with historical chemical references. More abstractly, Thomas conjectured in 1981 that the presence of a \emph{`negative loop'} in the system is \emph{necessary} for stable periodic behavior. In contrast, Ivanova \cite{Iv79} generalized \emph{positive-feedback cycles} into \emph{`critical fragments'} of the network and conjectured a complementary condition for oscillations: the existence of a critical fragment that does not involve all species in the network. See Mincheva and Roussel  \cite{MR07} for a review of Ivanova's work in English.  In the same paper, the authors however note that \emph{oscillations can
occur because of factors with different combinations of positive and negative cycles}. Gatermann et al. \cite{Gat2005} also acknowledge that oscillations can arise both from positive and negative feedback. In their own words, they \emph{`distinguish the type of instability according to the sign of the underlying
feedback loop. If the loop is positive the system is called autocatalytic. [...] An unstable chemical reaction with only a negative loop with at least three elements is called a nonautocatalytic oscillator.'}

In this section, we interpret our results from such historical perspective. We argue that the first of our criteria can be interpreted as the presence of unstable positive feedback within a stable subnetwork, while the second as the presence of unstable negative feedback. To keep this presentation simple, we proceed informally. We refer the interested reader to \cite{VasHunt, VasStad23}, where the following ideas have been discussed in full generality.

Again, let us start with pure linear algebra. We call \emph{unstable core} an invertible $n\times n$ matrix with a negative diagonal, which is unstable and such that none of its principal submatrices is unstable. {An $n\times n$ unstable core is further called an \emph{unstable-positive feedback} if its determinant has sign $(-1)^{n-1}$, and \emph{unstable-negative feedback} if it has sign $(-1)^n$.} As a direct consequence, unstable-positive feedbacks have one single real positive eigenvalue, while unstable-negative feedbacks have no real positive eigenvalues \cite[Lemma 6.3]{VasStad23}. Such definition generalizes positive and negative feedback cycles as, e.g.,
$$F^+:=\begin{pmatrix}
    -1 & 0 & 1 \\
    1 & -1 & 0 \\ 
    0 &  1 & -1 \\
\end{pmatrix} \quad \text{and} \quad F^-:=\begin{pmatrix}
    -1 & 0 & -1 \\
    1 & -1 & 0 \\ 
    0 &  1 & -1 \\
\end{pmatrix}, $$
where `positive' vs `negative' typically refers to the product of the off-diagonal entries. 

{
Assisted by the examples above, we proceed by making a few observations that, while limited and possibly imprecise, may help in understanding the design of chemical oscillators. The instability of positive feedback cycles, such as $F^+$, depends on a determinant calculation. Thus, this instability can be achieved by ensuring that the product of the off-diagonal entries is greater than that of the diagonal entries,} e.g.
$$F_u^+:=\begin{pmatrix}
    -1 & 0 & 2 \\
    1 & -1 & 0 \\ 
    0 &  1 & -1 \\
\end{pmatrix},\quad\text{with eigenvalues $(  -1.63 \pm 1.09i,0.26)$}.$$
The instability of such type of positive feedback cycle, in particular, does not depend on its length, i.e. on the size of the matrix. Consider now a stable $4\times4$ matrix with a negative diagonal, such that $F^+_u$ is one of its principal matrices, e.g.
\begin{equation}\label{eq:as}
A_s:=
\begin{pmatrix}
    -1 & 0 & 2 & 0\\
    1 & -1 & 0 & 2\\
    0 & 1 & -1 & 0\\
    0 & 0 & -1 & -1\\
\end{pmatrix},\quad\text{with eigenvalues $(-1,-1,-1,-1)$} .
\end{equation}
$A_s$ is indeed stable and can be seen as an overlap of the positive feedback cycle $F_u^+$ and the negative feedback cycle in $A[(2,3,4)]$. By construction, it is not a $P^-_0$ matrix because one of its principal submatrices is an unstable-positive feedback. A stable matrix that is not a $P^-_0$ matrix is precisely what our Criterion I requires $B(\bar{\mathbf{v}})$ to be. We underline the assonance with Ivanova's condition: if $A$ is not a $P^-_0$ matrix but is stable, then the unstable principal submatrix $A[\kappa]$ is necessarily a strict submatrix of $A$, i.e., $k<n$.

In striking contrast, instability for negative feedback cycles {may} depend on the length of the cycle. For instance, a negative feedback cycle with diagonal entries $F^-_{ii}=-1$ and off-diagonal entries $F^-_{i(i-1)}=1$ for $i=2,...,n$ and $F^-_{1n}=-2$ becomes unstable only for $n\ge 8$:
\begin{equation}\label{eq:f-u}F^-_u:=\begin{pmatrix}
      -1&     0&     0&     0    & 0&     0 &    0 &   -2\\
     1   & -1 &    0   &  0     &0   &  0   &  0 &    0\\
     0    & 1&    -1&     0    & 0    & 0   &  0&     0\\
     0     &0    & 1 &   -1    & 0     &0   &  0     &0\\
     0&     0   &  0  &   1   & -1     &0   &  0    & 0\\
     0 &    0  &   0   &  0  &   1    &-1   &  0   &  0\\
     0  &   0 &    0    & 0 &    0    & 1   & -1  &   0\\
     0   &  0&     0     &0&     0     &0   &  1 &   -1\\
\end{pmatrix},\end{equation}
with eigenvalues approximately $(-2.01 \pm 0.42i,-1.42 \pm 1.01i, -0.58 \pm 1.01i,
  0,0075\pm0.41i)$. We refer to such a heuristic observation as
$$\text{Negative-feedback instability {may require} length!}$$
Note also that $F^-_u$ is indeed an unstable $P^-$ matrix as Criterion II requires $B(\bar{\mathbf{v}})$ to be. 
We refer to \cite{MPSmith90, F19} for an extended analysis of the spectrum of feedback cycles. 

We now return to reaction networks and interpret such matrices as stoichiometric matrices of a (sub)network. The negative-diagonal condition simply reflects the fact that the species $X_i$ corresponding to the $i^{th}$ row is a reactant to the reaction $J(i)$ corresponding to the $i^{th}$ column. For a general algebraic treatment based on injective `Child-Selection maps' $J$, we refer again to \cite{VasHunt,VasStad23}. Negative-diagonal stoichiometric matrices have been called \emph{Child-Selection (CS-) matrices}. In its full generality, omitted here, such an approach does not depend on the chosen labeling of the network and naturally generalizes to CS-matrices that do not necessarily have a negative diagonal, depending on the presence of explicitly autocatalytic reactions. Moreover, \emph{stoichiometric autocatalysis}, as defined by Blokhuis et al. \cite{blokhuis20}, has been proven equivalent to the presence of an unstable-positive feedback that is a Metzler matrix, i.e., has nonnegative off-diagonal entries \cite[Thm. 7.3]{VasStad23}.
With a degree of consistency to Gatermann et al. \cite{Gat2005}, autocatalytic instabilities are then only due to the presence of a special type of unstable-positive feedback, while unstable-negative feedbacks always identify non-autocatalytic instabilities. 

We look at the stoichiometry of a network with these lenses: we aim to derive conclusions about the dynamical stability of a steady state. We put mass action systems aside for a moment, and we consider at first a more general Michaelis--Menten nonlinearity:
\begin{equation}\label{eq:MM}
    r_i=a_i\prod_{m=1}^{|S|}\bigg(\frac{x_m}{1+b_i^m x_m}\bigg)^{s^i_m},\quad a_i>0, b_i^m\ge 0.
\end{equation}
Mass action is recovered as the limit case when $b_i^m=0$ for any $m$ and $i$. The advantage of including also the $b_i^m$ parameters is that the mere presence of an unstable core implies the existence of a choice of parameters $a_i$, $b_i^m$ such that the system admits an unstable steady-state \cite[Corollary 5.1]{VasStad23}. Moreover, the presence of stoichiometry like $A_s$, \eqref{eq:as}, (which we refer to as \emph{unstable-positive feedback within a stable subnetwork}) and $F^-_u$, \eqref{eq:f-u}, (\emph{unstable-negative feedback}) is already a sufficient condition for the presence of a steady-state with purely imaginary eigenvalues of the Jacobian \cite[Corollary 5.17]{VasHunt}, pointing at periodic oscillations. The mathematical argument is analogous to Criteria I and II, respectively. The advantage of considering Michaelis--Menten over mass action is precisely that the sufficient conditions can be expressed at the level of naked stoichiometry, without taking into consideration the steady-state constraints imposed by the flux vector $\bar{\mathbf{v}}$. 

Such arguments identically apply to a general class of monotone chemical functions named \emph{parameter-rich} \cite{VasStad23}. Parameter-rich kinetics include e.g. Michaelis--Menten, Hill, and Generalized mass action kinetics. Unfortunately, mass action kinetics is not parameter-rich:
the restriction to  mass action, $b_i^m\equiv0$ for any $m$ and $i$, is then more challenging. However, it has been proved \cite{Ba-07,VasHunt} that the absence of unstable-positive feedbacks implies that the Jacobian is a $P^-_0$ matrix. Thus, Criterion I necessarily requires an unstable-positive feedback. Furthermore, \cite{VasStad23} also conjectured that the absence of $D$-unstable CS-matrices characterizes networks that only admit a unique and locally stable steady-state for any parameter choice. We may then still interpret Criteria I and II in terms of positive/negative feedback, but the presence of stoichiometries like $A_s$ and $F^-_u$ is no longer a sufficient condition for the assumptions of Criteria I and II to hold.

At present, the above observations offer thus only guidance on how to find periodic oscillations in mass action systems. In the next section, based on the above intuitions, we present three examples that admit periodic orbits.

\section{Examples}\label{sec:examples}
 The first example builds on an unstable-positive feedback within a stable subnetwork, i.e. stoichiometry of the type of $A_s$ \eqref{eq:as}, and applies Criterion I. The second and third examples are centered on unstable-negative feedbacks, i.e. stoichiometry of the type of $F^-_u$ \eqref{eq:f-u}, and apply Criterion II.
 \subsection{Example I: unstable-positive feedback within a stable subnetwork}
\begin{figure}[h!]
    \centering
    \includegraphics[width=1\linewidth]{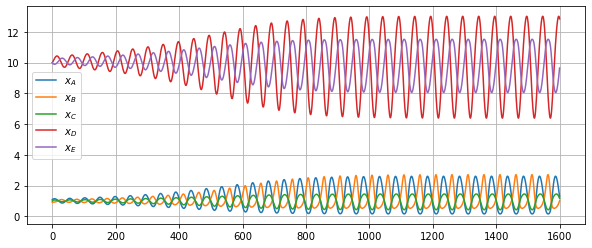}
    \caption{Numerical simulations for system \eqref{eq:ex1}. The steady-state flux vector has been chosen as $\bar{\mathbf{v}}=(1,1,1,1,1)$. The values of the (unique) unstable steady-state are  $\bar{x}=(1,1,1,10,10)$, which imply a chice of rates $(k_1,k_2,k_3,k_4,k_5)=(1,1,1,0.1,0.1)$. Initial conditions have been chosen $x(0)=(1.1,0.9,1,10,10)$. The plot shows convergence to a stable limit cycle.}
        \label{fig:1}
\end{figure}

Consider the following network with 5 species and 5 reactions:
\begin{equation}
\begin{cases}
A+B \quad &\underset{1}{\rightarrow} \quad C\\
B+C \quad &\underset{2}{\rightarrow} \quad E\\
C \quad &\underset{3}{\rightarrow} \quad A+D\\
D \quad &\underset{4}{\rightarrow} \quad 2B\\
E \quad &\underset{5}{\rightarrow} \quad C\\
\end{cases}
\end{equation}
The stoichiometric matrix and the kinetic matrix read:
\begin{equation}
N=\begin{pmatrix}
-1 & 0 & 1 & 0 & 0\\
-1 & -1 & 0 & 2 & 0\\
1 & -1 & -1 & 0 & 1\\
0 & 0 & 1 & -1 & 0\\
0 & 1 & 0 & 0 & -1    
\end{pmatrix}\quad Y=\begin{pmatrix}
1 & 0 & 0 & 0 & 0\\
1 & 1 & 0 & 0 & 0\\
0 & 1 & 1 & 0 & 0\\
0 & 0 & 0 & 1 & 0\\
0 & 0 & 0 & 0 & 1    
\end{pmatrix}
\end{equation}
Note the presence of an unstable-positive feedback
\begin{equation}\begin{pmatrix}\label{ex1:PF}
    -1 & 0 & 1 \\
-1 & -1 & 0 \\
1 & -1 & -1 
\end{pmatrix}, \quad\text{with eigenvalues approx $(0.32, 
  -1.66 \pm 0.56i),$}
\end{equation}
as a principal submatrix of a stable negative-diagonal matrix:
$$
\begin{pmatrix}
-1 & 0 & 1 & 0 \\
-1 & -1 & 0 & 2\\
1 & -1 & -1 & 0\\
0 & 0 & 1 & -1
\end{pmatrix},\quad\text{with eigenvalues approx $( -0.34 \pm 0.56i,
  -1, -2.32)$}.
$$
The unique steady-state flux vector is
$$\bar{\mathbf{v}}=(c,c,c,c,c), \text{ for $c\in \mathbb{R}_{>0}$,}$$
{
and the unique left kernel vector $w$ of the stoichiometric matrix $N$, 
$$w=(c,0,c,0,c) \text{ for $c\in \mathbb{R}_{>0}$,}$$
identifies the unique conserved quantity $\Omega$:
$$\Omega:=x_A+x_C+x_E.$$
Note also that the network is closed, i.e., it does not possess inflows or outflows.} The associated mass action system is:
\begin{equation}\label{eq:ex1}
\begin{cases}
\dot{x}_A=-k_1x_Ax_B+k_3x_C\\
\dot{x}_B=-k_1x_Ax_B -k_2x_Bx_C+2k_4x_D\\
\dot{x}_C=k_1x_Ax_B-k_2x_Bx_C-k_3x_C+k_5x_E\\
\dot{x}_D=k_3x_C-k_4x_D\\
\dot{x}_E=k_2x_Bx_C-k_5x_E\\
\end{cases}
\end{equation}
Without loss of generality, we fix $c=1$ and thus $\bar{\mathbf{v}}=(1,1,1,1,1)$ and $\operatorname{diag}(\bar{\mathbf{v}})=\operatorname{Id}$. We can then easily compute $B(\bar{\mathbf{v}})$:

\begin{equation}
\begin{split}
B(\bar{\mathbf{v}})=N \operatorname{Id} Y^T&=\begin{pmatrix}-1 & 0 & 1 & 0 & 0\\
-1 & -1 & 0 & 2 & 0\\
1 & -1 & -1 & 0 & 1\\
0 & 0 & 1 & -1 & 0\\
0 & 1 & 0 & 0 & -1    
\end{pmatrix}
\begin{pmatrix}
1 & 0 & 0 & 0 & 0\\
1 & 1 & 0 & 0 & 0\\
0 & 1 & 1 & 0 & 0\\
0 & 0 & 0 & 1 & 0\\
0 & 0 & 0 & 0 & 1 
\end{pmatrix}^T\\&=\begin{pmatrix}
-1 & -1 &  1 & 0  & 0\\
-1 & -2 & -1 & 2  & 0\\
1  & 0 & -2 &  0  & 1\\
0  &  0 &    1 &   -1 &    0\\
0  &   1&     1&     0&    -1
     \end{pmatrix},
     \end{split}
\end{equation}
with eigenvalues  \begin{equation}\label{eq:eigenvaluespF}(\lambda_1,\lambda_2,\lambda_3,\lambda_4,\lambda_5)\approx ( -3, -2.6, -1, -0.38, 0).
\end{equation}
Note the trivial eigenvalue zero due to the presence of one conserved quantity. On the other hand, in correspondence with the above analysis on the stoichiometric matrix, the  principal minor
\begin{equation}
    \begin{pmatrix}
        -1 & -1 & 1 & 0\\
        -1 & -2 & -1 &2\\
        1 & 0 & -2 & 0\\
        0 & 0& 1 & -1
    \end{pmatrix}
\end{equation} is stable, with eigenvalues $(\lambda_1,\lambda_2,\lambda_3,\lambda_4)\approx (-2.73 \pm  0.68i, -0.27 \pm 0.23i);$ and corresponding to the unstable positive-feedback \eqref{ex1:PF}, there is still a 3-principal minor with sign $(-1)^{3-1}=1$,
\begin{equation}\label{eq:pFex}
    \operatorname{det}\begin{pmatrix}
    -1 & -1 & 1\\
    -1 & -2 & -1\\
    1 & 0 & -2
\end{pmatrix}=1,
\end{equation}
and thus $B(\bar{\mathbf{v}})$ is not a $P^-_0$ matrix. Our Criterion I implies nonstationary periodic orbits. The equilibrium $\bar{x}=(1,1,1,1,1)$ is stable in its stoichiometric compatibility class, because of \eqref{eq:eigenvaluespF}. Intuitively, to make the instability of \eqref{eq:pFex} dominant, and thus achieve a Hopf bifurcation, we need to choose $h_4=1/\bar{x}_D$ and $h_5=1/\bar{x}_E$ small enough, i.e., $\bar{x}_D$ and $\bar{x}_E$ big enough. For example, the equilibrium $\bar{x}=(1,1,1,10,10)$ is unstable and encircled by a stable limit cycle. We confirm this via numerical simulations, see Figure \ref{fig:1}. 

\subsection{Example II and Example III: unstable-negative feedback}
\paragraph{Example II}
For the second case study, we elaborate on a family of examples presented by Claude Hyver in \cite{hyver1978}. 
The examples comprise $n+2$ species and $n+2$ reactions: 
\begin{equation}
\begin{cases}
\underset{F}{\rightarrow}\quad A_1 \quad \underset{1}{\rightarrow}\quad A_2 \quad \underset{2}{\rightarrow}\quad...\quad \underset{n-1}{\rightarrow}\quad A_n \quad \underset{n}{\rightarrow}\quad B+C\\
\quad \quad \quad\quad \quad \quad A_1+B \quad \underset{n+1}{\rightarrow}\quad ...\\
\quad \quad \quad\quad \quad \quad A_1+C  \quad \underset{n+2}{\rightarrow}\quad   ...
\end{cases}
\end{equation}
The actual example discussed in \cite{hyver1978} included a few more species, probably to enhance its chemical realism. We focus on the above simplification since it already contains all the mathematical features of our interest. The associated mass action system is:
\begin{equation}\label{eq:ex2}
\begin{cases}
\dot{x}_{A_1}=F-k_1x_{A_1}-k_{n+1}x_{A_1}x_{B}-k_{n+2}x_{A_1}x_{C}\\
\dot{x}_{A_2}=k_1x_{A_1}-k_2x_{A_2}\\
\vdots\\
\dot{x}_{A_n}=k_{n-1}x_{A_{n-1}}-k_{n}x_{A_{n}}\\
\dot{x}_B=k_{n}x_{A_{n}}-k_{n+1}x_{A_1}x_{B}\\
\dot{x}_C=k_{n}x_{A_{n}}-k_{n+2}x_{A_1}x_{C}
\end{cases}
\end{equation}

The stoichiometric matrix and the kinetic matrix are of the form:
\begin{equation}
    N=\begin{pmatrix}
        1 & -1 & 0 & ... & 0 & -1 & -1\\
        0 & 1 & -1 & ... & 0 & 0 & 0\\
        0 & 0 & 1 & ... & 0 & 0 & 0\\
         & & & \vdots & & &\\
        0 & 0 & 0& ... & -1 & 0 & 0\\
         0 & 0 & 0& ... & 1 & -1 & 0\\
         0 & 0 & 0& ... & 1 & 0 & -1\\
    \end{pmatrix}\quad\text{and}\quad
    Y=\begin{pmatrix}
        0 & 1 & 0 & ... & 0 & 1 & 1\\
        0 & 0 & 1 & ... & 0 & 0 & 0\\
        0 & 0 & 0 & ... & 0 & 0 & 0\\
         & & & \vdots & & &\\
        0 & 0 & 0& ... & 1 & 0 & 0\\
         0 & 0 & 0& ... & 0 & 1 & 0\\
         0 & 0 & 0& ... & 0 & 0 & 1\\
    \end{pmatrix}
\end{equation}
Without the first inflow column, the stoichiometric matrix is a $P^-$ matrix for any $n$, {and becomes unstable for $n\ge 7$. It is thus} an unstable-negative feedback for $n\ge 7$ large enough. Note that it is just a variant of the negative feedback cycles $F^-_u$, \eqref{eq:f-u}, where
the principal submatrix 
$$\begin{pmatrix}
    -1 & -2\\
    1 & -1 
\end{pmatrix},$$
in \eqref{eq:f-u} has been replaced by the principal submatrix:
$$\begin{pmatrix}
    -1 & -1 & -1\\
    1 & -1 & 0\\
    1 & 0 & -1
\end{pmatrix}.$$
 The unique steady-state flux vector is 
$$\bar{\mathbf{v}}=(3c,c,c,...,c,c,c).$$
Choosing $c=1$ without loss of generality, we can compute $B(\bar{\mathbf{v}})$:
\begin{equation}
\begin{split}
   B(\bar{\mathbf{v}})&= N\operatorname{diag}(3,1,1,...,1,1,1)Y^T\\
   &=\begin{pmatrix}
         -3 & 0 & ... & 0 & -1 & -1\\
         1 & -1 & ... & 0 & 0 & 0\\
         0 & 1 & ... & 0 & 0 & 0\\
         & & \vdots & & &\\
         0 & 0& ... & -1 & 0 & 0\\
          -1 & 0& ... & 1 & -1 & 0\\
          -1 & 0& ... & 1 & 0 & -1\\
         \end{pmatrix}.
   \end{split}
\end{equation}

$B(\bar{\mathbf{v}})$ is a $P^-$ matrix for any $n$. 
To confirm this, we argue by induction on $n\ge 2$, as it is enough to note that an increase in $n$ only adds a column and a row
$$\begin{pmatrix}
& & & 0 &&&\\
    &&&\vdots&&&\\
   0&...&0& -1&0&...&0\\
    &&&1&&&\\
   &&& \vdots&&&\\
    &&&0&&&
\end{pmatrix}.$$
Thus, it boils down only to check  whether
\begin{equation}
\begin{pmatrix}
-3 &0 & -1 & -1\\
1 & -1 & 0 & 0\\
-1 & 1  & -1 & 0 \\
-1 & 1 & 0 & -1\\
\end{pmatrix}
\end{equation}
is a $P^-$ matrix, which can be done by explicit computation. The case $n=1$ must be checked independently, as $B(\bar{\mathbf{v}})$ takes the different form:
$$B(\bar{\mathbf{v}})=\begin{pmatrix}
    -3 & -1 & -1\\
    0 & -1 & 0\\
    0 & 0 & -1
\end{pmatrix},$$
which is a $P^-$ matrix, since it is upper triangular with negative diagonal.

\begin{figure}
    \centering
    \includegraphics[width=1\linewidth]{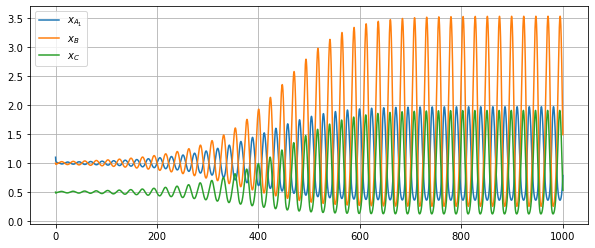}
    \caption{Numerical simulations for system \eqref{eq:ex2},with $n=10$. The steady-state flux vector has been chosen as $\bar{\mathbf{v}}=(3,1,1,...,1,1)$. The values of the (unique) unstable steady-states are $\bar{x}=(1,1,...,1,0.5)$, which imply a choice of rates $(F,k_1,...,k_{12})=(3,1,1,...,1,2)$. We have opted to choose $\bar{x}_B\neq\bar{x}_C$ because the two trajectories would fully overlap otherwise. Initial conditions have been chosen $x(0)=(1.1, 1,1,...,1,0.5)$. The plot shows convergence to a stable limit cycle. For graphical simplicity, we only include trajectories for $A_1,B,C$. }
    \label{fig:2}
\end{figure}

On the other hand, for $n$ large enough, the stability of $B(\bar{\mathbf{v}})$ is lost, and our Criterion II implies the existence of periodic orbits. We do not analytically prove here that $B(\bar{\mathbf{v}})$ is definitely unstable for $n$ large enough. We just check  explicitly that for $n<10$ $B(\bar{\mathbf{v}})$ is stable, while for $n=10,11,12,13$ a pair of complex conjugate eigenvalues with positive real part appears:
\begin{equation*}
\begin{cases}
n=1,...,9:\quad\text{all eigenvalues with negative real part;}\\
n=10:\quad \text{ one pair of eigenvalues with positive real part:}\quad 0.0049 \pm 0.2631i;\\
n=11:\quad \text{ one pair of eigenvalues with positive real part:}\quad 0.0091 \pm 0.2424i;\\
n=12:\quad \text{ one pair of eigenvalues with positive real part:}\quad  0.0119 \pm 0.2248i;\\
n=13:\quad \text{ one pair of eigenvalues with positive real part:}\quad  0.0139 \pm 0.2096i.\\
\end{cases}
\end{equation*}
We note a monotonicity in the values of the unstable eigenvalues. This family of networks confirms that the length of the structure helps periodic oscillations to appear. In Figure \ref{fig:2}, numerical simulations for $n=10$, and initial conditions nearby the (unstable) equilibrium $\bar{x}=(1,1,1,...,1,0.5),$ shows convergence to a stable limit cycle.

\paragraph{Example III} For the third case study, we generalize Example C from \cite{VasStad23} into a family of networks with $n+1$ species, $A_1,...,A_n,B$, $n>1$, and $2n+1$ reactions:
\begin{equation}
    \begin{cases}
      \underset{F_k}{\rightarrow}A_k,  \quad\quad\quad  2A_k+A_{k+1}\quad\underset{k}{\rightarrow}\quad2A_k+B\quad\quad\text{$k \in \mathbb{Z}^n$}\\
      B\underset{n+1}\rightarrow
      \end{cases}
\end{equation}
In \cite{VasStad23}, numerical simulations for $n=5$ showed the presence of a stable limit cycle. 
Note that species $A_1,...,A_n$ together with reactions $1,...,n$ give rise to a stoichiometry:
\begin{equation}\label{ex3:nf}
\begin{pmatrix}
0 & 0 & ... & 0 & -1\\
-1 &0 & ... & 0 & 0\\
0 & -1 &...& 0 & 0\\
&& \vdots &&\\
0 & 0 &... & -1 & 0
\end{pmatrix},
\end{equation}
which is indeed a $P^-_0$ matrix for any $n$ and unstable for $n\neq 2,4$: eigenvalues can be expressed in terms of roots of unit. The matrix does not have a negative-diagonal since all reactions are explicitly autocatalytic, but in the full generality addressed in \cite{VasStad23} such matrix is also an unstable-positive feedback for $n\neq 2,4$. Again, the length of the cycle helps instability to arise. The associated mass action system reads:
\begin{equation}\label{ex3:ma}
\begin{cases}
\dot{x}_{A_k}=F_k-a_{k-1}x_{A_{k-1}}^2x_{A_k}\quad\text{for $k\in \mathbb{Z}^n$},\\
\dot{x}_{B}=\sum_{k=1}^n a_{k-1}x_{A_{k-1}}^2x_{A_k} - a_{n+1}x_{B}.
\end{cases}
\end{equation}
Fixing all inflow constants $F_k\equiv 1$ yields a steady-state flux vector 
$$\bar{\mathbf{v}}=(1,1,1,...,1,1,n).$$
We compute 
\begin{equation}B(\bar{\mathbf{v}})=N\operatorname{diag}(\bar{\mathbf{v}})Y^T=
\begin{pmatrix}
-1 & 0 & ... & 0 & -2 & 0\\
-2 &-1 & ... & 0 & 0 & 0\\
0 & -2 &...& 0 & 0 &0\\
&& \vdots &&&0\\
0 & 0 &... & -2 & -1 &0\\
0 & 0 & ... & 0 & 0 &-n
\end{pmatrix},
\end{equation}
which is a $P^-$ matrix. $B(\bar{\mathbf{v}})$ also inherits the instability of \eqref{ex3:nf} for $n\ge5$. Criterion II implies nonstationary periodic orbits for $n\ge5$. Note that the form \eqref{ex3:ma} suggests the applicability of a Poincaré-Bendixson criterion for monotone cyclic feedback cycles \cite{MPSmith90}: we do not pursue here this direction.

\section{Fully-open systems}\label{sec:fullyopen}
%We conclude this paper with a result on fully-open networks.  % The reaction network $\mathbf{\Gamma'}=(S,R')$ obtained from $\mathbf{\Gamma}=(S,R)$ by considering an extended reaction set $R'\supseteq R$ with the addition of the missing inflow and outflow reactions is called \emph{fully-open extension} of $\mathbf{\Gamma}$. 
Albeit rather irrealistic, fully-open networks have been considered in relevant reaction network literature by Shinar \& Feinberg, \cite{ShiFei12} Banaji \& Craciun \cite{BaCra10}, Banaji \& Pantea \cite{BaPa16}, Conradi \& Pantea \cite{CoCa19}, among others. Our result characterizes the capacity of fully-open networks for Hopf bifurcation in terms of their capacity to admit a steady state with a complex-conjugate pair of eigenvalues with a positive-real part. For simplicity, we focus on  \emph{local Hopf bifurcation}. We recall the standard theorem \cite{GuHo84}.
\begin{theorem}[local Hopf bifurcation]\label{thm:Hopf}
Suppose that a parametric system
$$\dot{x}=f(x,\beta),\quad x\in\mathbb{R}^n,\;\beta\in\mathbb{R},$$
has a steady-state $(\bar{x}^*, \beta^*)$ at which the following two conditions hold:
\begin{enumerate}
\item the Jacobian $f_x(x)|_{(\bar{x}^*, \beta^*)}$ has a simple pair of purely imaginary eigenvalues $\lambda_{1,2}$ and $n-2$ eigenvalues with nonzero real part;
\item $\Re (\lambda_{1,2}(\beta^*))' \neq 0.$
\end{enumerate}
Then there is a local curve of steady-states $\bar{x}(\beta)$ that changes stability at $\bar{x}^*=\bar{x}(\beta^*)$. For some $\beta$ values in a neighborhood of $\beta^*$, the system admits periodic orbits.
\end{theorem}
The localization of the periodic orbits in the neighborhood of $\beta$ depends on the so-called \emph{Lyapunov coefficients}, which we do not discuss here. If $\dot{x}=f(x,\beta)$ satisfies the two conditions of Theorem \ref{thm:Hopf}, we say that the system undergoes a local Hopf bifurcation at $\beta=\beta^*$. We characterize local Hopf bifurcations for fully-open networks.

\begin{theorem}[Hopf bifurcation for fully-open networks]\label{thm:1}
Consider a fully-open reaction network $\mathbf{\Gamma}$ with associated mass action system \ref{maineq}. The following are equivalent:
\begin{enumerate}
    \item There is a choice of reaction rates {constants} such that the system undergoes a local Hopf bifurcation.
    \item There is a choice of reaction rates {constants} such that the system possesses a positive steady-state $\bar{x}$ such that its Jacobian matrix $f_x(x)|_{x=\bar{x}}$ has a simple pair of complex conjugate eigenvalues $\lambda_{1,2}$ with positive real part $\Re(\lambda_1)=\Re(\lambda_2)>0$ and no further eigenvalues $\lambda_i$, $i\neq 1,2$, with $\Re(\lambda_i)=\Re(\lambda_1)$. 
\end{enumerate}
\end{theorem}
\proof
The direction ($1\Rightarrow 2$) is trivial, as it is just one of the conclusions of Theorem \ref{thm:Hopf}. For the other direction ($2\Rightarrow 1$), consider the positive steady state $\bar{x}$ of $\mathbf{\Gamma}$ satisfying the assumptions. Introduce the positive bifurcation parameter $\beta\in\mathbb{R}_{>0}$ and define 
$$
\begin{cases}
\Delta(\beta):=\beta \operatorname{Id_{|S|}};\\
\mathbf{F}(\beta):=\Delta(\beta)\bar{x}.
\end{cases}
$$
Above, $\operatorname{Id_{|S|}}$ indicates the $|S|\times|S|$ identity matrix. We can then consider the following perturbation of $f(x)$
$$\dot{x}=h(x,\beta)=f(x)+\mathbf{F}(\beta)-\Delta(\beta)x.$$
Clearly, $h(x,0)=f(x)$. Moreover, $h(x,\beta)$ is fully-open as $f$: we interpret the vector $\mathbf{F}(\beta)$ as a perturbation of the inflow rates of $f$ and the diagonal entries of $\Delta(\beta)$ as a perturbation of the outflow rates of $f$. Moreover, by construction,
$$h(\bar{x},\beta)=f(\bar{x})+\mathbf{F}(\beta)-\Delta(\beta)\bar{x}=f(\bar{x})=0,$$
i.e., $\bar{x}$ is a steady state of $h(\beta)$ for any choice of $\beta$. For such perturbed choice of reaction rates, however, the Jacobian $h_x(x,\beta)$ of $h$ at $\bar{x}$ reads
$$h_x(x,\beta)|_{x=\bar{x}}=f_x(x)|_{x=\bar{x}}-\beta \operatorname{Id}_{|S|},$$
i.e., its spectrum corresponds to the spectrum of $f_x(x)|_{x=\bar{x}}$ translated to the left by $\beta$. Consider now the bifurcation value $\beta^*:=\Re(\lambda_1)>0$. By construction and assumption,  $h_x(x,\beta)|_{x=\bar{x}}(\bar{x},\beta^*)$ possesses a simple pair of purely imaginary eigenvalues $\mu_1, \mu_2$ and no other eigenvalues with zero real part. Moreover, $\mu_1(\beta), \mu_2(\beta)$ cross the imaginary axis transversely, i.e., $\Re(\mu_{1,2}'(\beta^*))'\neq 0.$ Indeed: simply note that $\mu _{1,2}(\beta)=\lambda_{1,2}-\beta$, and thus $|\mu_{1,2}'(\beta)| \equiv 1 \neq 0$. Thus, we can apply Theorem \ref{thm:Hopf}, and the fully-open system admits a local Hopf bifurcation. 
\endproof

%assume that $\bar{x}$ is a positive steady-state of the fully-open mass action system 
%$$\dot{x}=f(x),$$
%i.e. $f(\bar{x})=0$. We can consider now the following choice of parameters:

 %where $\mathbf{F} \in \mathbb{R}^{|S|}_{>0}$ is a constant vector and $D$ is  a $|S|\times|S|$ positive diagonal matrix $D>0$. Clearly, $h$ is of the same form as $f$: we can see $F$ as a perturbation of the inflow constant vector of $f$ and $D$ as a perturbation of the outflow rate of $f$. 

\begin{remark}For consistency with the rest of the paper, Theorem \ref{thm:1} is stated for mass action systems. A similar proof holds for any fully-open system where inflow rates are {free constant parameters and the steady-state derivative of the outflow rates can be chosen large ad libitum, e.g. if the outflow rates are linear:  a standard assumption in the dynamical modeling of biochemical systems}. In particular, Theorem \ref{thm:1} also holds if the outflow rates follow Generalized Mass Action, Michaelis--Menten, or Hill kinetics. The proof requires a minimal adaptation in the choice of rates to compute the Hopf bifurcation point. To do this, we refer again to the concept of parameter-rich kinetics \cite{VasStad23}.
\end{remark}

\section{Conclusion}\label{sec:conclusion}

We have stated two complementary criteria that guarantee the insurgence of nonstationary periodic orbits in mass action systems via a Hopf bifurcation. The criteria rely on the concept of $P^-$ matrices and are thus based on the sign of principal minors of the Jacobian. They consequently offer an evident computational advantage in comparison to the usual Hurwitz computation. More in detail, via Stoichiometric Network Analysis, we have expressed the Jacobian matrix $Jac$ evaluated at a steady state $\bar{x}$ as $$Jac=B(\bar{\mathbf{v}})\operatorname{diag}(1/\bar{x}_i),$$
where $\bar{\mathbf{v}}$ is a steady-state flux vector. The first criterion requires that $B(\bar{\mathbf{v}})$ is stable but not a $P^-_0$ matrix, while the second requires that $B(\bar{\mathbf{v}})$ is an unstable $P^-$ matrix. Moreover, we have interpreted the underlying chemical mechanisms as an unstable-positive feedback within a stable network (Criterion I) and an unstable-negative feedback (Criterion II). 

We have presented three examples where the criteria have been put into practice. Example I is a closed network with 5-species and a 3-species unstable-positive feedback: Criterion I proves the occurrence of a Hopf bifurcation as soon as the concentrations of the 3 species involved in the unstable-positive feedback become large in comparison to the others: this way the unstable-positive feedback becomes dominant and drives the dynamics towards an unstable region with the appearance of a stable limit cycle. Examples II and III are built around unstable-negative feedbacks of size $n$ and apply Criterion II. We have noted that the instability of negative feedback appears only for a sufficiently large size.  Consequently, both examples are families of networks of any size $n$: this confirms the validity of our criteria especially to address large networks. Periodic orbits appear for $n$ large enough: $n\ge 10$ in Example II, $n\ge 5$ in Example III. 

Finally, as an independent observation, we have proved that the capacity for Hopf bifurcation of fully-open systems is equivalent to the capacity for an unstable steady-state with a simple pair of eigenvalues with positive-real part.

{
To conclude the paper, we outlook three possible lines of research for the future. 
\paragraph{Purely stoichiometric criteria}
A direct inspection of the stoichiometric matrix \( N \) reveals that \emph{unstable-positive feedbacks within a stable subnetwork}, along with \emph{unstable-negative feedbacks}, may lead to purely imaginary eigenvalues of the Jacobian at the steady states of reaction systems with parameter-rich kinetics, such as Michaelis--Menten \cite{VasHunt,VasStad23}. The criteria established here for mass-action systems may provide similar insights but only apply to \( B(\bar{\mathbf{v}}) \) for a suitable choice of steady-state flux vector \( \bar{\mathbf{v}} \). The relationship between the bare stoichiometry \( N \) and the assumptions of these criteria requires further investigation. Ultimately, the goal is to identify straightforward stoichiometric patterns that sufficiently support the criteria, which would enhance our understanding of the chemical mechanisms driving oscillations in mass action systems.
\paragraph{Local stability analysis} In Section \ref{sec:linearalgebra}, we have argued that  standard local Hopf bifurcation (Theorem  \ref{thm:Hopf})
could be applied to establish periodic orbits in mass action systems following the same criteria exposed in this paper, upon positive answer of the following linear algebra question.}

{
\emph{$Q^*$: Let $A$ be an invertible matrix. Assume there exist two positive diagonal matrices $D_1,D_2$ such that  $\operatorname{inertia}AD_1\neq\operatorname{inertia}AD_2.$ 
Does it always exist a choice of path $\gamma(t)\in\mathbb{R}^n$ with $t\in[0,1]$ and $\operatorname{diag}\gamma(0)=D_1$, $\operatorname{diag}\gamma(1)=D_2$, and a value $t^* \in (0,1)$ such that
\begin{enumerate}
\item $A\operatorname{diag}(\gamma(t^*))$ has a \emph{simple} pair of purely imaginary eigenvalues $\lambda_{i,j}$;
\item $\Re(\lambda_{i,j}(t^*))'\neq 0$?
\end{enumerate}}}

{
In other words: it is always possible to choose the path connecting $D_1$, $D_2$ so that the standard local Hopf bifurcation theorem applies? This approach could assert better the regularity of the resulting periodic orbit and theoretically allow - albeit computationally expensive in practice - for the examination of its local stability via an analysis of the so-called \emph{first Lyapunov coefficient}.
\paragraph{Parameter regions for oscillations}
 On the other hand, a key advantage of using the theory of global Hopf bifurcation over local Hopf bifurcation is the availability of a topological `global' description of the continuum of periodic orbits, see the concept of \emph{snakes} by Mallet--Paret and Yorke \cite{MpYorke82} and the description therein. This framework could be in principle exploited to better understand parameter regions exhibiting periodic oscillations. However, establishing conclusions about the persistence of periodic behavior away from the bifurcation point requires a bound on the minimal period, which is generally a challenging task.\\
}

{\noindent \textbf{Acknowledgments:} I thank Balázs Boros, Bernold Fiedler, Vilmos Gáspár, Josef Hofbauer, Nidhi Kaihnsa, and János Tóth for useful and pleasant discussions. In particular,
János Tóth made me aware of Hyver’s Example II. This work has been supported by the
DFG (German Research Society), project n. 512355535.}

\bibliographystyle{plain}
\bibliography{references}
\end{document}